\documentclass[12pt]{article}
\usepackage{amsmath}
\usepackage{amssymb}
\usepackage{amsthm}

\def\sd{\partial}
\def\xyb{(\bar x,\bar y)}

\def\Pi{{\cal P}}

\def\N{\mathbb{N}}
\def\R{I\!\!R}
\def\la{\lambda}
\def\al{\alpha}
\def\del{\delta}
\def\ep{\varepsilon}
\def\sig
\def\sd{\partial}
\def\xb{\overline x}

\def\yb{\overline y}

\def\dom{{\rm dom}~}
\def\gr{{\rm Graph}~}

\def\rra{\rightrightarrows}

\def\ep{\varepsilon}
\def\la{\lambda}

\def\al{\alpha}

\def\vf{\varphi}
\def\del{\delta}

\def\Ga{\Gamma}
\def\ga{\gamma}

\def\sig{\sigma}

\def\cll{{\mathcal L}}

\def\cn{{\mathcal N}}

\def\cs{{\mathcal S}}

\def\ck{{\mathcal K}}
\def\cm{{\mathcal M}}

\def\dom{{\rm dom}~}

\def\Sur{{\rm Sur}}
\def\sur{{\rm sur}}
\def\reg{{\rm reg}}

\def\cl{{\rm cl}}

\def\dfn{\noindent{\bf Definition}~}

\def\thm{\noindent{\bf Theorem}~}
\def\crl{\noindent{\bf Corollary}~}
\def\rem{\noindent{\bf Remark}~}
\def\prop{\noindent{\bf Proposition}~}
\newcommand{\black}{\hspace{10pt}\rule{6pt}{6pt}}
\begin{document}

\title{\huge A Sard Theorem for  Tame Set-Valued Mappings
}

\author{\large\bf A. Ioffe\thanks{Department of Mathematics, Technion,
Haifa 32000, Israel.} }
\date{}

\maketitle

\begin{abstract} If $F$ is a set-valued mapping from $\R^n$ into
$\R^m$ with closed graph, then $y\in \R^m$ is a critical value of
$F$ if for some $x$ with $y\in F(x)$, $F$ is not metrically
regular at $(x,y)$. We prove that the set of critical values of a
set-valued mapping  whose graph is a definable (tame) set in an
$o$-minimal structure containing additions and multiplications is
a set of dimension not greater than $m-1$ (resp. a porous set). As
a corollary of this result we get that the collection of
asymptotically critical values of a semialgebraic set-valued
mapping has dimension not greater than $m-1$, thus extending to
such mappings a corresponding  result by Kurdyka-Orro-Simon for
$C^1$ semialgebraic mappings. We also give an independent proof of
the fact  that a definable continuous real-valued function is
constant on components of the set of its subdifferentiably
critical points, thus extending to all definable functions a
recent result of Bolte-Daniilidis-Lewis for globally subanalytic
functions.
\end{abstract}

\section{Introduction.}
The classical Sard (or Morse-Sard) theorem states that the
collection of critical values of a $C^k$-mapping $F$ from (an open
subset of) $\R^n$ into $\R^m$ has Lebesgue measure zero, provided
$k\ge\max\{n-m+1,1\}$. The fundamental role of the Sard theorem in
analysis and differential geometry comes from the fact that for a
regular (non-critical) value $y$ of $F$ the set of solutions of
the equation $F(x)=y$ (if nonempty) is a nice set (a manifold)
which responds to variations of the right-hand side in a stable
and non-chaotic way. The Sard theorem therefore ascertains that a
typical value of a sufficiently smooth mapping is regular.

Such a result would be highly welcome in variational analysis in
which the stability issue is of an extreme importance. It would be
highly desirable to be able to make similar statements e.g. about
systems of inequalities or other relations of interests in
variational analysis.

At the first glance this does not seem to be possible. The Sard
theorem is sharp and there are widely known examples (Whitney
\cite{HW}, Yomdin \cite{YY3}) showing that for a less smooth
function or mapping Sard's theorem does not hold. The most precise
result was proved by Bates \cite{SB}: Sard's theorem holds for
$C^{n-m,1}$-mappings ($n-m$ times continuously differentiable with
locally Lipschitz derivatives of order $n-m$). Here ``Lipschitz''
cannot be stregthened to ``H\"older'' as was found by Norton
\cite{AN}.

However, recently Kurdyka-Orro-Simon \cite{KOS} proved that the
collection of critical and asymptotically critical values of a
semialgebraic $C^1$-mappings is a semialgebraic set of dimension
$m-1$ or less. Several results of Morse-Sard-type were proved for
real-valued functions under even more general assumptions:
\cite{AR,YY1} (quantitative results for maximum and minimax of
smooth families of functions), \cite{LR} (distance function to a
$C^{\infty}$-submanifold of a Riemann manifold), \cite{DAC}
(generalized critical values of a $C^1$-function definable in an
o-minimal structure), \cite{BDL,BDL1,BDLS} (critical points of
globally subanalytic functions). The last three papers have
largely stimulated this study.

These results demonstrate that the differentiability requirement
can be substantially weakened in exchange for some structural
restrictions on the class of mappings or functions. The main
result of this paper  shows that the frameworks of this trade-off
can be considerably expanded. This is the statement of the main
theorem.

\vskip 2mm

\thm 1. {\it If $F:\ \R^n\rightrightarrows \R^m$ is a tame
set-valued mapping with locally closed graph, then the set of
critical values of $F$ is a $\sigma$-porous\footnote{A set $Q$ in
a metric space is called {\it porous} if there is a $\la
>0$ such that for any $x\in Q$ and any $r>0$ the set
$B(x,r)\backslash Q$ contains a ball of radius $\la r$. A
$\sigma$-porous set is a countable union of porous sets. A
$\sigma$-porous set in $\R^n$ is both of the first Baire category
and Lebesgue measure zero.} set in $\R^m$. In particular it  has
Lebesgue measure zero. Moreover, if the graph of $F$ is a
definable set, then the set of critical values of $F$ is also a
definable set of dimension not exceeding $m-1$.}

\vskip 1mm

The next two sections contain all information from variational
analysis and the theory of o-minimal structures which is necessary
for the proof of the theorem. Here we shall only add a couple of
general remarks.

First we note that the concept of a `critical value'' provided by
modern variational analysis is very natural. (Actually this
concept seems to be defined here for the first time. However the
``parent'' concept of (metric) regularity has been thoroughly
studied during last two decades, see e.g. \cite{AI00,BM,RW}.) We
observe that restricted to single-valued continuously
differentiable mappings, this definition reduces to the classical
concept: $y$ is a critical value of $F$ if there is an $x$ such
that $F(x)=y$ and the rank of the derivative $\nabla F(x)$ is
smaller than $m$.

Definability and tameness are fundamental concepts of the theory
of o-minimal structures (see e.g. \cite{MC,LVD,BT}) which is also
being very actively developed last two decades, partly in response
to Grothendieck's call for a new ``tame topology'' based on ''a
real look ... at the context at which we live, breathe and work''.
\cite{AG}.  The important point about definable and tame objects
is that they are void of ``pathologies" so typical for generic
objects of nonsmooth analysis (e.g. Lipschitz functions that
cannot be recovered from their subdifferentials).

A search for good classes of nonsmooth functions (``less subject
to wildness'', I would add, again quoting  \cite{AG}) for which
certain results could be proved, was among the dominant themes in
nonsmooth analysis since practically its very beginning. Just
mention lower $C^2$ functions \cite{RT}, amenable functions,
prox-regular functions \cite{RW}, semismooth functions \cite{RM},
minimal cuscos \cite{BZ}, partially smooth functions \cite{AL}.
Functions of these classes well serve for the purposes they have
been created, but none of the classes have structural properties
compared to the properties of definable and tame functions and
sets (which need not be differentiable or even continuous). Thanks
to these properties, definable and tame functions and sets look
like an almost ideal playground for applicable finite dimensional
variational analysis.

 As a consequence of the
theorem, we recover the part of the theorem of Kurdyka, Orro and
Simon \cite{KOS} relating to  the dimension of the set of
asymptotic critical values of semialgebraic mappings (not the
fibration part of the theorem). Actually we get an extension of
this theorem to set-valued mappings with semialgebraic graphs. We
get this result as a part of a more general theorem in which  a
``stratification'' of asymptotically critical values by rates of
asymptotic decline of the rates of surjection is taken into
consideration. We also give a separate proof of an $o$-minimal
extension of the recent result by Bolte, Danniilidis and Lewis
\cite{BDL} saying that a continuous globally analytic function
which is continuous on its domain is constant on connected
components of the set of its critical points.

Theorem 1 is proved in \S5. In \S4 we state and prove some
preliminary results needed for the proof of the theorem. Some of
them are probably new but some (e.g. definability of derivative)
are known. We give short proofs of the latter as well, just for
convenience. The principal results here are Proposition 1 (which
builds a bridge between the two parts by showing that the ``rate
of surjection'' which is a quantitative measure of regularity is a
definable or tame function, provided the graph of the set-valued
mapping is respectively definable or tame) and Proposition 6
(showing that for certain definable families of functions, uniform
smallness of functions implies smallness of derivatives on big
sets).   The last \S6 contains statements and proofs of the
corollaries mentioned  above.

\vskip 1mm

{\bf Acknowledgement}. This paper was written while I was on
sabbatical in the Department of Computer Science of Dalhousie
University. I wish to express my gratitude to the department and
especially to Jon Borwein for hospitality and excellent working
conditions I was provided with. I am also thankful to Adrian Lewis
for inspiring discussions.

\section{Openness, regularity and critical values.}
 The concept defined below
makes sense in every Banach, and actually in every metric space.
This level of generality is not needed here, so we define
everything for finite dimensional Euclidean spaces and refer the
reader to \cite{AI00} for the general theory.

So let again \ $F:\ \R^n\rightrightarrows \R^m$ \ be a set-valued
mapping. For $(x,y)\in\gr F$ we set
$$
{\rm Sur}F(x,y)(\la) = \sup\{ r\ge 0:\ y+ rB_Y\subset F(x+\la
B_X)\}.
$$
and then for $\xyb\in\cl(\gr F)$ (the closure of $\gr F$) define
the {\it rate of openness (surjection)} of $F$ at $\xyb$ by
$$
{\rm sur}F(\xb|\yb) = \liminf_{(x,y,\la)\to (\xb,\yb,+0)}
\frac{1}{\la}{\rm Sur}F(x,y)(\la),
$$
For single valued $F$ we usually write $\sur F(\xb)$ (instead of
$\sur F(\xb|F(\xb))$).

The function $\la\to \Sur F(x|y)(\la)$  is called the {\it modulus
of surjection} of $F$ at $\xyb$. Here, as usual, $B_X$ etc. is the
unit ball in $\R^n$ and we set $\sup \emptyset=0$ (or else, we can
calculate the liminf in the definition of ${\rm Sur} F$ only along
sequences of $(x,y)\in \gr F$).

It follows from the definition that the function $\sur F$ is
defined on the closure of $\gr F$. If however $F$ is a set-valued
mapping with closed graph, then  it is possible to show that
\footnote{This fact probably has not been explicitly mentioned
earlier but it easily follows from the slope characterization of
the rate of surjection given in \cite{AI00}}
$$
\sur F(\xb|\yb) = \liminf_{(x,y)\underset{{\rm Graph}~F}\to
(\xb,\yb)}\liminf_{\la\to+0}\frac{1}{\la}\Sur(x,y)(\la)
$$
In general, the quantity in the right-hand side of the equality
can by greater. Consider for instance the following mapping
$\R\rra\R$:
$$
F(x)= \{ y:\; 0<|y|<|x|\}.
$$
Then $\sur F(0,0)=0$ (take $x_n=\la_n=n^{-1}$, $y_n=n^{-2}$) but
the right-hand side quantity is equal to $\infty$.

The reciprocal of \sur$F(\xb|\yb)$ is the {\it rate of metric
regularity} of $F$ at $\xyb$.
$$
\reg F(\xb|\yb) = [\sur F(\xb|\yb)]^{-1}.
$$
This is a quantitative measure of stability of solution $\xb$  of
$y\in F(x)$ at $y=\yb$ for it is precisely the lower bound of
positive $K$ such that
$$
d(x,F^{-1}(y))\le K d(y,F(x))
$$
for all $(x,y)$ of a neighborhood of $\xyb$.

It is said that $F$ is {\it regular} at $(x,y)$ or that $(x,y)$ is
a {\it regular point} of $F$   if $\sur F(x|y)>0$.

Otherwise $(x,y)$ is  a {\it singular point} of $F$. Finally, $y$
is a  {\it singular} or {\it critical} value of $F$ if there is a
$x$ such that $y\in F(x)$ and $\sur F(x|y)=0$.

\vskip 1mm

\rem. Observe that the above definition of a critical value covers
both the case of a ``proper'' critical value when $(x,y)$ belongs
to the graph of $F$ and of a  ``generalized'' critical value when
$(x,y)$  belongs to the closure of $\gr F$ but not  to the
graph of $F$ itself. In principle, if we impose no topological
conditions on $F$, it may happen that the a critical value of $F$
(even proper) is a regular value of the set-valued mapping whose
graph is the closure of $\gr F$. Consider, for instance the
following set-valued mapping $F:\ \R^2\rra\R^2$:
$$
F(0)= \{ 0\}; \quad F(x)=\| x\|B\backslash ({\rm the}\  x\ {\rm
axis}),\; {\rm if}\;  x\neq 0.
$$
Then zero is a proper critical value of $F$ but a regular value
of the mapping whose graph is the closure of $\gr F$.

To avoid pathologies like that, we shall mainly consider
set-valued mappings with locally closed graphs. Another reason for
introducing such an assumption is that the all known regularity
criteria need it.

The following known facts will be used in the sequel (see
\cite{AI00})

\vskip 1mm

$\bullet$ \ $y$ is a critical value of $F$ if and only if it is a
critical value of the projection $\R^n\times\R^m\to \R^m$
restricted to $\gr F$;

\vskip 1mm

$\bullet$ \ $\sur F$ is a lower semicontinuous function;

\vskip 1mm

$\bullet$ \ if $F(x)=Ax$ is a linear  operator, then for any $x$
$$
\sur A(x) = \Sur A(x)(1) = \inf_{\| y^*\| = 1} \| A^*y^*\| = \|
A^{*-1}\|^{-1},
$$
so that we can write just $\sur A$ etc.;

\vskip 1mm

$\bullet$ \  if $F$ is single-valued and continuously
differentiable at $x$, then $\sur F(x) = \sur (\nabla F(x))$,
where by $\nabla F(x)$ we denote the Jacobian matrix of $F$ or/and
the corresponding linear operator $\R^n\to\R^m$.

\vskip 1mm

$\bullet$ \ if $F(x)=H(x) + A(x)$, where $A$ is a linear operator,
$H$ is a set-valued mapping with locally closed graph and $y\in
H(x)$, then $\sur F(x|y+A(x))\ge \sur H(x) - \| A\|$;

\vskip 1mm

$\bullet$ \ if $F=H\circ G $, where $G$ is a $C^1$ mapping into
the domain space of $H$ and the graph of $H$ is locally closed,
then
$$
\sur G(x)\cdot\sur H(G(x)|y)\le\sur F(x|y)\le \|\nabla
G(x)\|\cdot\sur H(G(x)|y).
$$

\vskip 1mm

The most general regularity criterion (actually, the precise
formula for the rate of regularity)  is based on the concept of
slope introduced by  DeGiorgi-Marino-Tosques. We shall state it
only for single valued mappings as it is sufficient here.   Let
$f$ be an extended-real-valued function which is finite at $x$.
The {\it slope} of $f$ at $x$ is the lower bound of $K\ge 0$ such
that $f(x)\ge f(u)-K\| x-u\|$ for all $u$ of a neighborhood of
$x$. (The usual convention $\inf\emptyset =\infty$ applies.) The
slope is usually denoted $|\nabla f|(x)$ to emphasize that for a
Fr\'echet differentiable function the slope coincides with the
norm of the gradient.

Given a mapping $F:\ \R^n\rra\R^m$ defined and continuous in a
neighborhood of a certain $\xb$, we set for any $y\in\R^m$
$f_y(x,v)= \| y-F(x)\|$.  Then {\it $\sur F\xyb$ is the upper
bound of $\ga\ge 0$ having the property that there is an $\ep>0$
such that $|\nabla f_y|(x)\ge \ga$ for all $x$ and $y$ satisfying
 $\| x-\xb\|<\ep$, $y\neq F(\xb)$}. (The same result will be
 obtained if, instead of all $y\neq F(\xb)$ we shall take $y\in
 U\backslash F(\xb)$, where $U$ is an arbitrary neighborhood of
 $F(\xb)$.)

Let us define the {\it slope of} $F$ at $x$ by
$$
{\rm Sl}~F(x) = \inf_{y\neq F(x)}|\nabla f_y|(x,v).
$$
Then the above stated result can be equivalently expressed as
$$
\sur F(\xb)=\liminf_{x\to \xb}{\rm Sl}~F(x).
$$
It follows in particular that
$$
\sur F(x)\le {\rm Sl}~F(x), \;\forall\; x\in\dom F.
$$

\section{$O$-minimal structures and definable functions.}
We give below the statements of main definitions and facts without
proofs. There are two excellent introductions to the subject:
\cite{MC,LVD}. General $o$-minimal structures can be associated
with various linearly ordered sets but we here shall consider only
structures associated with the real line $\R$. We shall keep the
notation  $\R^n$ also for  Euclidean spaces and denote the inner
product by $(\cdot|\cdot)$.

\dfn 1. An {\it  structure} on $\R$ is a sequence $\cs = (\cs_n)$,
$(n\in \N)$ such that for each $n$

\vskip 2mm

(D1) $\cs_n$ is a Boolean algebra of subsets of $\R^n$, that is,
$\emptyset\in\cs_n$ and $\cs_n$ contains unions, intersections and
complements of its elements;

\vskip 1mm

(D2) If $A\in\cs_n$, then $A\times\R$ and $\R\times A$ belong to
$\cs_{n+1}$;

\vskip 1mm

(D3) $\{ x=(x_1,...x_n): \; x_i=x_j\}\in\cs_n$ for any $1\le
i<j\le n$;

\vskip 1mm

(D4) If $A\in\cs_{n+1}$ then the projection $\pi(A)$ of $A$ to
$\R^n$ ($\pi:\; (x_1,...,x_n,x_{n+1})\to (x_1,...,x_n)$) belongs
to $\cs_n$;

\vskip 1mm

A structure is called {\it $o$-minimal} (short for ``order
minimal'') if in addition

 (D5) $\{(x,y)\in\R^2:\ x<y\}\in\cs_2$.

\vskip 1mm

 (D6) the elements of $\cs_1$ are precisely finite unions of
points and open intervals.

\vskip 1mm

The elements of $\cs_n, n=1,2,...$ are called {\it definable} (in
$\cs$). A set $Q$ is called {\it tame} if its intersection with
any bounded definable set is a definable set. A (set-valued)
mapping $F$ from a subset of $\R^n$ into $\R^m$ is called {\it
definable} ({\it tame}) if its graph is a definable (tame)   set.
Likewise, a real-valued function defined on a subset of $\R^n$ is
definable if its graph is a definable set in $\R^{n+1}$.

In variational analysis it is often convenient to work with
extended-real-valued functions defined on all of $\R^n$. The
definition can easily be extended to such functions: $f$ is
definable if its graph $\{(x,\al)\in\R^n\times\R:\ \al=f(x)\}$ is
definable along with the sets $\{ x:\ f(x)=\infty\}$ and $\{ x:\
f(x) = -\infty\}$.

\vskip 1mm

The subtle points of the definition are (D4) and (D6). (D4) is
usually the most difficult part of the proof that a certain
collection of sets is a structure. A consequence of it is that any
set obtained from definable sets with the help of finitely many
existential and universal quantifiers $\exists$, $\forall$
(applied to variables only, not to sets, functions and other
parameters) and boolean operations is also definable.

(D6) is basically  responsible for a number of remarkable
structural (tameness) properties of definable sets and functions
(e.g. Monotonicity, Cell Decomposition and Definable Choice
theorems stated below) which exclude any possibility of ``wild''
behavior.

Here are several examples of $o$-minimal structures most useful
for analysis.

1) Let us call a set $Q\subset \R^n$ an {\it open polyhedron} if
it is the intersection of finitely many open half spaces $\{ x: \
f(x)<0\}$ and hyperplanes  $\{x:\ f(x)= 0\}$, where $f(x) =
(a|x)+\al$ are affine functions. The structure $\cs_{\rm lin}$ of
{\it semilinear} sets is formed by finite unions of open
polyhedrons. All axioms of $o$-minimal structures are easy to
verify in this case.

\vskip 1mm

2) If we replace affine functions by  polynomials, we obtain the
structure $\cs_{\rm alg}$ of {\it semialgebraic sets}. Here again
all axioms are verified easily, with the exception of (D4). The
latter is the subject of a deep Tarski-Seidenberg theorem (see
\cite{BCR}). A consequence of this fact is that semialgebraic sets
admit {\it elimination of quantifiers}, that is any set obtained
from semialgebraic sets with the help of quantifiers and boolean
operations can also be obtained by means of a quantifier free
formula involving only level and sublevel sets of polynomials and
boolean operations.

\vskip 1mm

3) The above scheme no longer works if we make a step further and
replace polynomials by real analytic function. We can define {\it
semianalytic} sets in the same way but using arbitrary real
analytic functions. However, in this case (D4) and (D6) do not
hold. Indeed, the set of zeros of $\sin x$ is an infinite
collection of isolated points and there are also examples of
bounded semianalytic sets whose projections are not semianalytic.

Nonetheless there exists a rich $o$-minimal structure  in which
all bounded semianalytic sets are definable. A  set $Q\subset R^n$
is called ${\it subanalytic}$ if locally near each of its points
it is a projection of a bounded semianalytic set, that is if  for
any $x\in Q$ there is an open neighborhood $U$ of $x$ and a
bounded semianalytic $S\subset \R^{n+k}$ such that the projection
of $S$ to $\R^n$ coincides with $Q\bigcap U$. A set $Q\subset
\R^n$ is called {\it globally subanalytic} if $F(Q)$ is
subanalytic whenever $F$ is a semialgebraic homeomorphism of
$\R^n$ onto $(-1,1)^n$.

 It turns out that globally subanalytic sets already satisfy all
 axioms and the corresponding $o$-minimal structure is denoted
 $\cs_{\rm an}$.

 It is possible to give an alternative ``non-constructive''
 description for $\cs_{\rm lin}$ and $\cs_{\rm alg}$, namely $\cs_{\rm lin}$
 is the minimal structure containing graphs of affine functions
 and $\cs_{\rm alg}$ is the minimal structure containing graphs of all
 polynomials. It turns out that $\cs_{\rm an}$ is the minimal
 structure containing semialgebraic sets and graphs of restrictions of real analytic
 functions to balls.

\vskip 1mm

 4) The minimal structure $\cs_{\rm an,exp}$ containing all globally
 subanalytic sets and the graph of the exponent $e^x$ is also o-minimal.
Clearly
$$
\cs_{\rm lin}\subset\cs_{\rm alg}\subset\cs_{\rm
an}\subset\cs_{\rm an,exp}
$$

 \vskip 1mm

We shall further consider only $o$-minimal structures satisfying
the additional property

\vskip 1mm

 (D7) The graphs of addition: $\{(x,y,z)\in\R^3: \
z=x+y\}$ and multiplication $\{ (x,y,z)\in \R^3:\ z=x\cdot y\}$
belong to $\cs_3$.

The semi-linear structure does not belong to this class but every
semilinear set is semialgebraic, so all results valid for the
latter are also valid for semilinear sets.

The following are some simple properties  of definable and tame
sets and functions which are obtained from the axioms with
relative easiness.

$\bullet$ \ the closure and the interior of a definable (tame) set
is a definable (tame) set;

$\bullet$ \ a function is definable (tame) if and only if its
epigraph $\{ (x,\la):\ \la\ge f(x)\}$ and hypograph $\{ (x,\al):\
a\le f(x)\}$ are definable (tame) sets;

$\bullet$  \ the derivative (also partial) of a definable (tame)
function is a definable (tame) function;

$\bullet$ \ under a suitable agreement about operations with
infinite values (e.g. $\infty -\infty=\infty; 0\cdot\infty=\infty$
etc. ) the collection of definable (tame) extended-real-valued
functions is stable under summation, subtraction, multiplication
and operations of pointwise maximum and minimum; composition of
definable mappings is definable.

 $\bullet$ \ If
$f_1,...f_k$ are definable (tame) functions and $\R^n$ is
partitioned into $k$ definable sets $X_1,...,X_k$, then the
function $f$ equal to $f_i$ on $X_i$ is definable(tame).

$\bullet$ \ the functions
$$
\inf_y f(x,y)\quad {\rm and}\quad  \sup_y g(x,y)
$$
are definable, provided so are $f, \ g$;

$\bullet$ \ the image and the preimage of a definable set under a
definable mapping is a definable set; the image of a tame set
under a proper tame mapping is a tame set.

\vskip 2mm

And now several fundamental results characterizing the tameness
properties of definable sets and functions.

\vskip 2mm

 \noindent {\bf Monotonicity theorem}. {\it Let $f$ be a
definable function on $\R$. Then $\dom f$ is a finite union of
points and (open) intervals, and on each of the intervals $f$ is
either constant or strictly monotone and continuous.}

\vskip 2mm

\noindent {\bf Uniform finiteness theorem}. {\it Let $F$ be a
definable set-valued mapping from $\R^n$ into $\R^m$. Suppose that
every $F(x)$ contains finitely many points. Then there is a
natural $N$ such that the number of points in every $F(x)$ does
not exceed $N$.}

The next theorem uses the concept of a {\it cell} whose definition
we omit. Although this concept will be often used in what follows,
we do not need the specific structure of cells described in the
formal definition. For us it will be sufficient to think of a
$C^k$-{\it cell} of dimension $r$ as of an $r$-dimensional
$C^k$-manifold which is the image of the cube $(0,1)^r$ under a
definable $C^k$-diffeomorfism. As follows from the definition, an
$m$-dimensional cell in $\R^m$ is an open set.

\vskip 2mm

\noindent{\bf Cell Decomposition Theorem} \ {\it (a) Let $Q\subset
\R^m$ be a definable set. Then for any $k$,  $Q$ can be
represented as a disjoint union of a finite number of cells of
class $C^k$;

(b) Let $F$ be a definable mapping from a set $Q\subset \R^n$ into
$\R^m$. Then there exist a partition of $Q$ into a finite number
of cells of class $C^k$ such that the restriction of $F$ on each
cell is a mapping of class $C^k$. }

\vskip 1mm

The maximal dimension of the cell in a decomposition of a
definable set is called the {\it dimension} of the set. Of course,
the dimension does not depend on the choice of decomposition. The
important fact concerning dimension is that {\it the dimension of
the image of a definable set under a definable (single-valued)
mapping cannot be greater than the dimension of the preimage}.

Another consequence of the cell decomposition theorem is that any
definable set has a finite number of connected components. More
careful analysis leads to the conclusion that {\it any connected
definable set is definably pathwise connected}, that is any two
points of the set can be joined by a definable continuous curve
lying completely in the set.

For a set-valued mapping $F: \R^n\rightrightarrows \R^m$ we denote
$\dom F= \{x\in \R^n:\;\exists y [y\in F(x)]\}$, that is to say,
the projection of the graph  of $F$ onto the domain space. A {\it
selection} of $F$ is, as usual, a mapping $\vf(x)$ from $\dom F$
into the image space such that $\vf(x)\in F(x)$ for all $x\in\dom
F$.

\vskip 2mm

\noindent{\bf Definable Choice Theorem}. {\it Any definable (resp.
tame) set-valued mapping has a definable (resp. tame) selection}.

\vskip 1mm

We conclude the introduction with the following simple example
which demonstrates the difference between definability and
tameness and shows that one must be more careful when working with
tame objects.

The  function $\sin t$ is a tame function and  $t^{-1}$ is a
semialgebraic function but the composition $\sin t^{-1}$ is not
even a tame function (its restriction to e.g. $(0,1)$ is not
definable in any $o$-minimal structure). Thus a composition of
tame functions may be not a tame function. Observe that $\sin t$
is not a proper map.

\section{Some preliminary results.}

\prop 1. {\it Let $F$ be a definable (tame) set-valued mapping
from $\R^n$ into $\R^m$. Then $\sur F$ is a definable (resp. tame)
function with $\dom (\sur F)\subset \cl(\gr F)$.}

\proof We can represent the graph of $(x,y,\la)\to \Sur
F(x|y)(\la)$ (considered a function on $\gr F\times(0,\infty)$) as
$(P\bigcup Q)\bigcap S$, where

\vskip 2mm

\noindent$ P = \gr F\times(0,\infty)\times\{ 0\}$;
$$
Q\ =\ \{(x,y,\la,r):\; \forall \ 0\le\rho<r, \ \forall\ v\ [\|
v-y\|\le\rho\Rightarrow \ \exists u,\ \| u-x\|\le \la, v\in
F(u)]\};
$$
and
$$
S=\{(x,y,\la,r):\; \forall\ \ep>0\ \exists \ v \ [\|v-y\|<r+\ep \
\& \  \forall u\ [\| u-x\|\le\la\ \Rightarrow v\not\in F(u)]]\}.
$$
If $F$ is a definable mapping, then  $P$ is clearly definable and
 the other two sets are also definable by (D4) as was explained in
 the previous section.  Hence so is the
graph of $\Sur F(x|y)(\la)$. If $F$ is a tame mapping, then for
any $K>0$ the intersections of three sets with $\{ (x,y,\la,r):\;
\| x\|\le K,\ \| y\|\le K,\ 0\le \la,r\le K\}$ is a definable set,
so the intersection of the graph of $\Sur F$ with any such set is
a definable set.

We note further  the epigraph of $\sur F$ is  the intersection of
the closure of the epigraph of $\la^{-1}\Sur (x|y)(\la)$ with the
set $\{x,y,\la,\al):\ \la = 0\}$, so $\sur F$ is a definable (or
tame) function and its domain lies in the closure of $\gr F$ by
definition.\quad\black

\vskip 2mm

\prop 2. {\it Let  $F$ be a continuous definable (single-valued)
mapping from an open definable subset of $\R^n$ into $\R^n$. Then
the dimension of the set of critical values of $F$ is not greater
than $n-1$}.

\proof  By the cell decomposition theorem, $\dom F$ is the union
of $C^1$-cells and the restriction of $F$ onto each of them is
$C^1$. Then the set of critical values of the restriction of $F$
to any cell of dimension $n$ has measure zero   by the Sard
theorem, hence, by definability, its dimension cannot be greater
than $n-1$. On the other hand the image of the union of $F$-images
of all other cells is also a set of dimension not greater than
$n-1$. \qquad\qquad\qquad \qquad\qquad\qquad\qquad\qquad\qquad
\qquad\qquad\qquad\black

\vskip 2mm

\prop 3. {\it Let $U$ be an open definable subset of $\R^n$ and
$F$ a continuous definable (single-valued) mapping from $U$ into
$\R^m$. Assume that $\sur F(x)=0$ for every $u\in U$. Then $\dim
F(U)\le m-1$.}

\proof Assume the contrary: $\dim F(U)=m$. Then there is an
$m$-dimensional  cell $Q\subset F(U)$ which by definition is an
open subset of $\R^m$.

Applying the definable choice theorem, we shall find a definable
mapping $G: \ Q\to \R^n$ such that $\{(x,y):\; x=G(y),\ y\in
Q\}\subset \gr F$. This means that $F\circ G= Id|_Q$.

As $G$ is definable, there is a smaller definable cell $\tilde Q
\subset Q$ such that $G$ is $C^1$ on $\tilde Q$. We have $(F\circ
G)(\tilde Q)= \tilde Q$. Therefore the set $G(\tilde Q)$ also has
dimension $m$. This set is also definable as a definable image of
a definable set. Therefore there is an $m$-dimensional cell
$D\subset G(\tilde Q)$ such that the restriction of $F$ to $D$ is
$C^1$. But then $G\circ F|_D= Id|_D$ and as $F$ continuously
differentiable at every $x\in D$ and $G$ is continuously
differentiable at every $y= F(x), x\in D$. Without loss of
generality we can identify $D$ with $(-1,1)^m$. So we have $\nabla
G(y)\circ\nabla (F|_D)(x)=I$, that is $\nabla (F|_D)(x)$ is a
non-singular operator $\R^m\to \R^m$ which by the Lusternik-Graves
theorem means that $\sur (F|_D)(x)>0$. As $\sur F(x)\ge \sur
(F|_D)$ (because $D\subset \tilde Q$ and $F$ is $C^1$ on $\tilde
Q$), we get a contradiction.
\qquad\qquad\qquad\qquad\qquad\qquad\black

\vskip 1mm

\rem.  Observe that a continuously differentiable mapping into
$\R^m,\ m\ge 2$ with the rate of surjection identically zero can
be surjective (see \cite{YY2}).

\vskip 2mm

\prop 4 (Differentiability theorem). {\it If $f$ is a definable
function on $\R$, then $f'$ is also definable and

\vskip 1mm

(a) $\dom f'$ is an open set;

(b) $f'$ is continuous on $\dom f$;

(c) $\dom f\backslash\dom f'$ is a finite set}.

\proof We have
$$
\begin{array}{l}
\gr f'=\{(t,\al)\in\dom f\times\R:\ \forall\ \ep>0\ \exists \
\del>0\ \\  \qquad\qquad\qquad\qquad\qquad\qquad [|\tau-t|<\del
\Rightarrow \ |f(\tau)-f(t)-\al (\tau-t)|<\ep|\tau-t|]\}.
\end{array}
$$
Hence $\gr f'$ is a definable set.

(a) Let $t\in\dom f'$. Then for any $\del >0$ the sets
$(t,t+\del)\cap \dom f$ and $(t-\del,t)\cap\dom f$ are nonempty.
The monotonicity theorem now implies (as every monotone function
is almost everywhere differentiable) that also the sets
$(t,t+\del)\cap \dom f'$ and $(t-\del,t)\cap\dom f'$ are nonempty
for any positive $\del$, hence they are infinite. As these sets
are definable, by (D6) there must be a $\del>0$ such that
$(t,t+\del)\subset\dom f'$ and $(t-\del,t)\subset\dom f'$, that is
$(t-\del,t+\del)\subset\dom f'$.

(b) As $f'$ is a definable function, it obeys the monotonicity
theorem. Thus we only need to observe that if there are
$\tau<\tau'<\tau''$ such that $(\tau,\tau'')\subset\dom f'$ and
$f'$ is monotone on $(\tau,\tau')$ and $(\tau',\tau'')$, then
$$
\lim_{t\to\tau'-0}f'(t)=\lim_{t\to\tau'+0}f'(t)
$$
for otherwise $f$ would not be differentiable at $\tau'$.

(c) The set $\dom f\backslash\dom f'$ must have Lebesgue measure
zero due to almost everywhere differentiability of monotone
functions. Therefore it cannot contain intervals and, being
definable, must be finite. \qquad\qquad\black

\vskip 2mm

\crl {\it Let $f$ be a definable function on an interval $(a,b)$.
Then there is a finite number of points $a=t_0<t_1<...<t_k=b$ such
that on each interval $(t_i,t_{i+1})$ $f$ is continuously
differentiable and $f'$ either strictly positive, or strictly
negative or identically equal to zero.}

\proof Let $t_i$ be either the points of non-differentiability or
$f$ or isolated points of $\{t: \ f'(t)=0\}$ or the ends of
intervals  in the decomposition of the set according to (D6).

\vskip 2mm

\prop 5 (Uniform integrability lemma). \ {\it Let $Q\subset\R^n$
be a definable set, and let $\vf(x,t,s)$ be a definable function
on $Q\times(a,b)\times (0,1)$, where $-\infty<a<b<\infty$. Suppose
that $\vf(x,t,s)\to 0$ as $s\to 0$ uniformly on $Q\times (a,b)$.
Then}
$$
\int_a^b|\vf_t(x,t,s)|dt\to 0,\quad {\rm as}\; s\to\infty,\;
\quad{\rm uniformly\ on}\quad Q.
$$
Here $\vf_t$ is the derivative of $\vf$ with respect to $t$.

\proof By the preceding corollary for any $x\in Q$, $s\in (0,1)$
there are finitely many, say $N(x,s)$ points $\tau_i=\tau_i(x,s)$
(with $\tau_0=a,\ \tau_{N(x,s)}=b$) on $(a,b)$ such that on each
interval bounded by a pair of adjacent points, $\vf$ is
continuously differentiable with respect to $t$ and the derivative
is either identical zero or does not change the sign. The
set-valued mapping which associates with every $(x,s)$ the
collection of these points $\tau_i(x,s)$ is definable.

Indeed the set
$$
\{ (x,t,s)\in Q\times(a,b)\times (0,1):\; \vf(x,\cdot,s)\ {\rm is
\ discontinuous \  at}\; t\}
$$
is definable. Indeed, this set is equal to the intersection of
$$
\{(x,t,s):\; \exists r>0, \ \ep>0\ \forall \del\in(0,\ep)\
[\vf(x,t+\del,s)-r>\vf(x,t-del,s)]\}
$$
and
$$
\{(x,t,s):\; \exists r>0, \ \ep>0\ \forall \del\in(0,\ep)\
[\vf(x,t+\del,s)+r<\vf(x,t-del,s)]\}.
$$
The same is true for the points at which the derivative of $\vf$
with respect to $t$ is discontinuous. Equally simple arguments
lead to the conclusion that the set of $(x,t,s)$ such that
$\vf_t(x,t,s)=0$ (the derivative with respect to $t$) and $t$ is
either an isolated point of the zero set of $\vf_t(x,\cdot,s)$ or
an end of the interval at which the function is equal to zero, is
also definable.

By the Uniform Finiteness Theorem the numbers $N(x,s)$ are
uniformly bounded, that is there is an $N$ such that $N\ge N(x,s)$
for all $x\in Q,\ s\in (0,1)$.

Fix an $\ep>0$ and choose a $\del> 0$ such that
$|\vf(x,t,s|\le\ep/2N$ if $0<s<\del$. Therefore for any $x\in Q,\
s\in (0,\del)$
$$
\int_{\tau_i}^{\tau_{i+1}}|\vf_t(x,t,s)|dt =
|\vf(x,\tau_{i+1},s)-\vf(x,\tau_i,s)|\le\ep/N
$$
and therefore,
$$
\int_a^b|\vf_t(x,t,s)|dt\le\ep.\qquad\qquad\qquad\black
$$

We shall use this result as a starting point for obtaining uniform
estimates for the norms of derivatives of  ``small'' definable
mappings. Let $F$ be a $C^1$-mapping from a neighborhood of $x\in
\R^n$ into $\R^m$. We shall denote by $\nabla F(x)$ the Jacobian
matrix of $F$ at $x$, namely
$$
(\nabla F(x))_j^i= \frac{\sd F^i}{\sd x_j}(x),
$$
where $x=(x_1,...,x_n)$ and $F(x)= (F^1(x),...,F^m(x))$. Recall
also that the norm of a linear operator defined by matrix
$A=(a_j^i)$ in the standard basis of $\R^n$ is
$$
\| A\| = \Big(\sum_{i=1}^m\sum_{j=1}^n (a_j^i)^2\Big)^{1/2}.
$$

\vskip 1mm

\prop 6. {\it Let $Q= (-1,1)^n$  and let $F:\ Q\times (0,1)\to
\R^m$ be a definable mapping having the property that $\|
F(x,s)\|\to 0$ uniformly on $Q$ when $s\to 0$. Then for any
$\ep>0$ there is a $\del>0$ such that, whenever $s<\del$, there is
an open set $\Omega=\Omega(s)\subset Q$  such that $\|\nabla_x
F(x,s)||<\ep$
for all $x\in\Omega$}.\\
Here the subscript $x$ refers to differentiation with respect to
$x$ only, not to $s$.

\proof Set for simplicity $ \vf_{i}(x,s)=(F(x,s))^i$. By the
assumption every $\vf_i$ goes to zero uniformly in $x\in Q$ as
$s\to 0$. Set $Q_j= \{ (x_1,...,x_{j-1},x_{j+1},...,x_n):\
|x_k|<1\}$. Applying Proposition 5 for each $\vf_i$ and treating
consecutively every $x_j$ as $t$, we conclude that for every $i,j$
$$
\int_{-1}^1\Big |\frac{\sd\vf_i(x,s)}{\sd x_j}\Big | dx_j \to
0\quad {\rm uniformly\ on}\; Q_j.
$$
As an consequence we get
$$
\int_{Q}\|\nabla_x F(x,s)\|dx \to 0,\qquad{\rm as}\; s\to 0.
$$

Indeed by the cell decomposition theorem we can partition $Q$ into
a finite number of cells such that $F$ is continuously
differentiable on each of them. Clearly, continuous
differentiability on a cell of dimension $n$ is the same as
continuous differentiability without any restrictions. Equally
clear is the union of all cells of dimension $n$ is a set of full
measure in $Q$. This means that the integral above makes sense. On
the other hand
$$
\int_{Q}\|\nabla_x
F(x,s)\|dx\le\sum_{ij}\int_{Q_j}\Big(\int_{-1}^1 \Big
|\frac{\sd\vf_i(x,s)}{\sd x_j}\Big | dx_j\Big)
dx_1,...,dx_{j-1}dx_{j+1}...dx_n.
$$

 We can choose $\del>0$ so small that the above integral is
not greater than $\ep/2$ if $s<\del$. But then for any such $s$
the Lebesgue measure of $\{x\in Q: \ \|\nabla_x F(x,s)\|<\ep\}$
must be at least $1/2$. This is a definable set, hence it has
nonempty interior as its measure is
positive.\qquad\qquad\qquad\qquad\qquad\qquad\qquad\qquad\qquad\qquad\black

\section{Proof of Theorem 1.}
 \noindent{\bf Step 1}. There is no loss of generality in
assuming that the graph of $F$ is closed (so that any critical
value is ``proper''). Indeed, as the graph of $F$ is locally
closed, the closure operation does not add points to  a small
neighborhood of any point of the graph, so that any critical value
of $F$, no matter ``proper'' or ``generalized'' remains a critical
value of the mapping whose grapf is $\cl(\gr F)$.

Furthermore, as we mentioned in \S 3, $y$ is a critical value of
$F$ if and only if it is a critical value of the restriction of
the projection $\R^n\times \R^m\to\R^m$ to the closure of the
graph of $F$. It follows that it is sufficient to prove the
theorem for a single valued mapping $F$ which is a restriction of
a linear operator $A\in\cll(\R^n,\R^m)$ to a tame set $\dom
F\subset \R^n$.

We also observe that any cell in $\R^m$ of dimension smaller than
$m$ is a porous set (which is immediate from definitions) , hence
any definable set in $\R^m$ of dimension smaller than $m$ is
$\sigma$-porous. It follows that the theorem will be proved if we
show that {\it for any definable mapping which is a restriction of
a linear operator to a definable set the dimension of the set of
its critical values cannot exceed $m-1$}.

Indeed, the mapping $F$ restricted to the set $\{ x\in\dom F:\ \|
x\|\le N\}$ is definable by definition and any critical value of
$F$ is a critical value of the restriction if $N$ is sufficiently
big.

\vskip 5mm

\noindent{\bf Step 2}.   By Proposition 2 the theorem is true if
the dimension of $\dom F$ coincides with the dimension of the
image space. Assume now that for a given $m\ge 1$ the theorem
holds for any definable mapping whose domain has dimension not
greater than $r\ge m$ and let $F$ be a mapping from $\R^n$ into
$\R^m$ with $\dim(\dom F)=r+1$ \ (of course $n\ge r+1$).

By the cell decomposition theorem $\dom F$ can be partitioned into
finitely many $C^k$-cells $C_i$ \ ($k\ge r+2-m$). Denote by

$\cn=\{\bigcup C_i:\ \dim C_i= r+1\}$ \ the union of all
$(r+1)$-dimensional cells of the partition:

$\ck=\{\bigcup \big(\cl(C_i)\bigcap(\dom F)\big),\ \dim C_i\le
r\}$ \ the union of intersections of the closures of cells of
dimensions $\le r$ with the domain of $F$.

Then

(a) the collection of critical values of $F|_{\ck}=A|_{\ck}$ is a
set of dimension $\le m-1$ by the induction assumption;

(b) for any cell $C_i$  the collection of critical values of
$F|_{C_i}= A|_{C_i}$ is a set of dimension $\le m-1$ by the
classical Sard theorem.  Indeed, if $\dim C_i= s$, then $C_i$ is
the image of the open $s$-cube $(-1,1)^s$ under a $C^k$ mapping
$G$. Therefore singular points of $A|_{C_i}$ are among singular
points of $A\circ G$. The latter is a $C^k$ mapping from the cube
into $\R^m$ and as $k\ge s+1-m$ as $s\le r+1$, Sard's theorem
applies.

Let  $x\in\dom F$ be a singular point of $F$ which does not belong
to any of the two above mentioned types. This means that

\vskip 2mm

(c) \ $x$ is a regular point of the restriction of $F$ to the cell
of the partition containing $x$. We claim that for some $i$ there
is a sequence $(x^{\nu})\subset C_i$ converging to $x$ such that
$$
\lim_{\nu\to+0}\sur F|_{C_i}(x^{\nu}) = 0. \leqno (1)
$$
Indeed, as $\sur F(x)=0$,  there is a sequence $(x^{\nu})$
converging to $x$ and such that ${\rm Sl}~F(x^{\nu})\to 0$ and
there is no loss of generality in assuming that all $x^{\nu}$
belong to the same cell, call it $C$. Furthermore, it follows from
the definition of the slope and the inequality at the end of
Section 2 that
$$
{\rm Sl}~F(x^{\nu}) \ge {\rm Sl}~F|_C(x^{\nu})\ge\sur
F|_C(x^{\nu})
$$
which implies (1)

We note further that $x^{\nu}$ cannot belong to $\ck$ for in this
case  $x$ also belongs to $\ck$ as the latter is closed. Thus $x$
is a singular point of $F|_{\ck}$ which is the case of (a). Thus
$x^{\nu}\in\cn$ for all $\nu$.  The limiting point $x$ cannot
belong to the same cell as $x^{\nu}$ since  we assume that $x$ is
a regular point of the restriction of $F$ to the cell.

Thus, we arrive to the following situation : there is a cell $C$
of dimension $r+1$ or higher and a sequence $(x^{\nu})$ converging
to $x$ such that $x^{\nu}\in C$ for all $\nu$, \ $x\in (\dom
F)\backslash C$ and (1) holds.

Let $\cm$ denote the collection of such points $x$ (associated
with the same cell $C$). We have to show that

$$
  \dim F(\cm)\le m-1. \leqno (2)
$$

\vskip 3mm

\noindent{\bf Step 3}. Thus we reduce the problem to the
following. Given

 $\bullet$\ a definable mapping $F$ from $\R^n$ into $\R^m$
 which is a restriction of a linear operator $A$ to a definable set
 $\dom F\subset \R^n$;

$\bullet$ an  cell $C\subset \dom F$  of class $C^1$ (and
dimension $\ge r+1)$;

$\bullet$  a nonempty  set
$$
\cm=\{ x\in(\dom F)\backslash C:\;  \forall s>0\ \exists\ y\in C \
[\| x-y\|<s,\; \sur F(y)< s]\}.
$$
It is clear from the definition of $\cm$ that it is a definable
set lying completely in the closure of $C$, that is in the
boundary of $C$, as $\cm$ and $C$ do not meet.

We have to prove that (2) holds. Assume by way of contradiction
that $\dim A(\cm)=m$. Then $\cm$ contains a $q$-dimensional ($q\ge
m$) cell of class $C^1$ whose $A$-image has dimension $m$. This
means that there is a diffeomorphism $G$ of $Q=(-1,1)^q$ into
$\cm$ such that   $\dim (A\circ G)(Q)=m$. If $\sur (A\circ G)(u) =
0$ for all $u\in Q$ then by Proposition 3 $\dim (A\circ G)(Q)\le
m-1$, so the rate of surjection of $A\circ G$ must be positive at
certain points of $Q$. Since $\sur (A\circ G)$ is a lower
semicontinuous function, we can assume, taking a smaller cube if
necessary,  that
$$
\sur(A\circ G)(u)\ge\al> 0, \quad \forall\; u\in Q. \leqno (3)
$$
and also that $G$ satisfies the Lipschitz condition in $Q$.

\vskip 5mm

\noindent{\bf Step 4}. \ Consider the set
$$
\Ga=\{ (x,y,s):\; x\in \cm,\ y\in C,\ s\in (0,1),\; \| y-x\|<s, \
\sur F(y)<s\}.
$$

This is a definable set and its projection onto the $x$-component
space is $\cm$. By the definable choice theorem there is a
definable mapping $y(x,s)$ from $\cm\times (0,1)$ into $\R^n$ such
that $(x,y(x,s),s)\in \Ga$ for all $(x,s)$. We have: $\|
x-y(x,s)\|<s$ for all $x\in \cm$ and all $s\in (0,1)$. Set
$$
\Psi (x,s) = x-y(x,s);\quad \Phi(u,s)= \Psi (G(u),s).
$$
 Then $\Phi(u,s)\to 0$ uniformly on $Q$ as $s\to 0$.

Applying  Proposition 6 to $\Phi$, we conclude that for each given
$\ep>0$ there is an $s=s(\ep)\le\ep$ and an open set
$\Omega(\ep)\subset Q$ such that $\|\nabla_u \Phi(u,s)\|<\ep$ for
all $u\in \Omega (\ep)$. This means that
$$
\|\nabla G(u)- \nabla (y_{\ep}\circ G)(u)\|<\ep, \leqno (4)
$$
where $y_{\ep}(x)= y(x,s(\ep))$.

Now for any $\ep>0$ choose $u=u(\ep)\in \Omega(\ep)$ such that
$y_{\ep}$ be differentiable at $u$ . Then, the inequality above
implies that
$$
\nabla (A\circ G)(u)= \nabla (A\circ y_{\ep}\circ G)(u)) +T(\ep),
\leqno (5)
$$
where $T(\ep)=(A\circ(\nabla G -\nabla(y_{\ep}\circ G)))(u)$ is a
linear operator from $\R^q$ into $\R^m$ with $\| T(\ep)\|\to 0$ as
$\ep\to 0$. We have by (3)
$$
\sur(\nabla (A\circ G)(u)) =\sur(\nabla(F\circ G)(u))=\sur(F\circ
G)(u)\ge \alpha. \leqno (6)
$$
On the other hand, by (4)
$$
\begin{array}{lcl}
\sur(\nabla(A\circ y_{\ep}\circ G))(u)&=&\sur(\nabla (F\circ
(y_{\ep}\circ G))(u)\\ &=& \sur(F\circ (y_{\ep}\circ G))(u)\le
(K+\ep)\sur F(y_{\ep}(G(u)),
\end{array}
$$
where $K= \|\nabla G(u)\|$.

The last three relations are however contradictory, as $\sur
F(y_{\ep}(G(u)))\to 0$ and from  (5), (6) we get
$$
0<\al\le \sur(\nabla (A\circ G)(u))\le\sur (\nabla (A\circ
y_{\ep}\circ G)(u))) +\| T(\ep)\|\to 0.
$$
This completes the proof of the theorem.

\section{Some corollaries.}

In \cite{KOS} Kurdyka, Orro and Simon proved that the set of
asymptotically critical values of a continuously differentiable
semialgebraic mapping $\R^n\to\R^m$ has dimension less  than $m$.
Their proof  was based on  calculation of an  estimate for  the
$m$-dimensional measure of the set.  Theorem 1 allows to avoid
these calculations and to get some information of asymptotical
critical values of  definable set-valued mappings in general.
First we note the following

\vskip 2mm

\prop 7. {\it Let $H: \ \R^n\rra \R^m$, and let $\rho (x)$ be
continuously differentiable strictly positive function. Set
$L(x)=H(\rho(x)x)$. Then for any $(x,y)\in\gr L$
$$
\sur L(x|y)\le \big(\rho(x)+\| \rho'(x)\|\cdot\| x\|\big)\sur
H(\rho(x)x|y).
$$
\proof Indeed, the norm of the derivative of the mapping $x\to
\rho (x)x$ at $x$ is not greater than $\rho (x)+ \|\rho'(x)\|
x\|$.
\qquad\qquad\qquad\qquad\qquad\qquad\qquad\qquad\qquad\qquad\black

\vskip 1mm

Let $\eta(t)$ be a strictly positive continuous function on
$[0,\infty)$ such that
$$
\int_0^{\infty}\frac{1}{\eta(t)}<\infty. \leqno (7)
$$

Given a set valued mapping $F: \ \R^n\rra \R^m$, we call $y\in
\R^m$ an {\it asymptotically $\eta$-critical value} of $F$ if
there is a sequence of pairs $(x^{\nu},y^{\nu})$ such that
$y^{\nu}\in F(x^{\nu})$,\ $\| x^{\nu}\|\to\infty$, \ $y^{\nu}\to
y$  and $\eta(\|x^{\nu}\|)\cdot\sur F(x^{\nu}|y^{\nu})\to 0$. We
shall denote by $K_{\infty}(F,\eta)$ the set of asymptotically
$\eta$-critical values of $F$.

\vskip 2mm

\thm 2. {\it Let $F$ be a definable set-valued mapping with
locally closed graph, and let $\vf(t)$ be strictly increasing
continuously differentiable positive definable function on
$[0,\infty)$, bounded from above and equal to zero at $0$. Set
$$
\eta(t) = \frac{1}{\vf'(t)}
$$
Then $K_{\infty}(F,\eta)$ is a  definable set with
$\dim(K_{\infty}(F,\eta))<m$.}

\proof Clearly, $\eta$ satisfies (7). We obviously have
$$
\begin{array}{l}
K_{\infty}(F)=\{y\in\R^m:\ \forall N\in(0,\infty),\ \ \exists
(u,v)\in \R^{n}\times\R^m
\\\qquad\qquad\qquad[v\in
F(u),\ \| v-y\|<N^{-1},\ \| u\|\ge N, \ \eta(\| u\|)\cdot\sur
F(u|v)<N^{-1}]\}
\end{array}
$$
which shows that the set is definable by Proposition 1.

Without loss of generality we may assume that $\vf(t)\to 1$ as
$t\to\infty$. Let $\psi (\cdot)$ be the inverse of $\vf(\cdot)$,
that is $\psi(\vf((t))\equiv t$. Set
$$
G(u) = F(\psi(\| u\|)\frac{u}{\| u\|}).
$$
In other words, we consider the following pair of mutually inverse
changes of variables:
$$
x= \psi(\|u\|)\frac{u}{\| u\|};\quad u = \vf(\| x\|)\frac{x}{\|
x\|}
$$
which transfer the open unit ball into the entire space and vice
versa. We obviously have for each pair of corresponding $x$ and
$u$:
$$
\| x\| =\psi (\| u\|),\quad \| u\| = \vf(\| x\|).
$$
(Here and below we consider the Euclidean norm in $\R^n$.)

We can apply Proposition 7 to get estimates of the rate of
surjection of $G$: for each $u\neq 0$ with $\| u\|<1$ and each
$y\in G(u)$
$$
\sur G(u|y)\le \big(\frac{\psi(\| u\|)}{\| u\|} +
\big\|\big[\frac{\psi(\| u\|)}{\| u\|}\big]'\big\|\cdot\|
u\|\big)\sur F((x|y).
$$
After simple calculation we get:
$$
\big\|\big[\frac{\psi(\| u\|)}{\| u\|}\big]'\big\|\cdot\| u\|\le
\|\psi'(\| u\|)+ \frac{\psi(\| u\|)}{\| u\|},
$$
so that
$$
\sur G(u|y)\le\ 2\big(\|\psi'(\| u\|)+ \frac{\psi(\| u\|)}{\|
u\|}\big)\sur F(x|y)= 2(\frac{1}{\vf'(\| x\|)} + \frac{\|
x\|}{\vf(\| x\|)}\big)\sur F(x|y).
$$

We now recall that $\eta(\cdot)$ is reciprocal of $\vf'(\cdot)$
and that by (7) $\eta(t)$ grows to infinity faster than $t$. Thus,
we can be sure that for sufficiently large $x$
$$
\sur G(u|y)\le 3\eta(\| x\|)\cdot\sur F(x|y). \leqno (8)
$$

Let now $\tilde G$ be set-valued mapping whose graph is the
closure of $\gr G$. Clearly, if $\| x^{\nu}\|\to\infty$ then the
sequence of the corresponding $u^{\nu}$ contains a subsequence
converging to an element of the unit sphere. Therefore as follows
from (8) every asymptotical $\eta$-critical value of $F$ is a
critical value of $\tilde G$. The latter is a definable set-valued
mapping, so by Theorem 1 the entire set of its critical values,
including $K_{\infty}(F,\eta)$ has the dimension smaller than $m$.
\qquad\qquad\qquad\qquad\qquad\qquad\qquad\qquad\qquad\black

\vskip 1mm

Let us call, following \cite{KOS},  a point $y\in \R^m$ an {\it
asymptotically critical value} of $F$  if there is a sequence of
pairs $(x^{\nu},y^{\nu})$ such that $y^{\nu}\in F(x^{\nu})$,\ $\|
x^{\nu}\|\to\infty$, \ $y^{\nu}\to y$  and $\|x^{\nu}\|\sur
F(x^{\nu}|y^{\nu})\to 0$. We shall denote by $K_{\infty}(F)$ the
set of asymptotically critical values of $F$.

Of course, any asymptotically critical value is asymptotically
$\eta$-critical for any $\eta(t)$ satisfying the requirements in
the definition but we cannot, in principle, expect an
asymptotically $\eta$-critical value to belong to $K_{\infty}(F)$.
However, Kurdyka, Otto and Simon showed in \cite{KOS} (Lemma 3.1)
that in case when $F$ is a semialgebraic $C^1$-mapping on $\R^n$
(or an open semialgebraic subset of $\R^n$), there is an $\ga>0$
(depending on $F$) such that $K_{\infty}(F,\eta)\subset
K_{\infty}(F)$ if $\eta(t)= t^{1+\al}$ with $0<\al<\ga$.

The proof of this fact given in \cite{KOS} extends without change
to arbitrary semialgebraic set-valued mappings. Actually, the only
use of differentiability in the proof is in the definition of the
distance from the derivative of $F$ to the set of singular
operators playing the same role in the definition of
asymptotically critical value given in \cite{KOS} as the rate of
surjection in the definition above. But this distance is precisely
the rate of surjection of $F$ at the corresponding point -- the
fact well known in variational analysis and actually also proved
in \cite{KOS} (Propositions 2.1, 2.2).

Thus, combining the quoted result of \cite{KOS} with Theorem 2, we
get an extension of (the first part of) the main theorem of
\cite{KOS} to arbitrary semialgebraic set-valued mappings.

\vskip 2mm

\thm 3. {\it Let $F$ be a semialgebraic set-valued mapping from
$\R^n$ into $\R^m$ whose graph is locally closed . Then
$K_{\infty}(F)$ is a closed semialgebraic set of dimension smaller
than $m$.}

\vskip 1mm

Our final result is an extension to  definable functions of a
recent theorem of Bolte-Daniilidis-Lewis \cite{BDL} stating that a
continuous globally subanalytic function is constant on every
connected component of its ``subdifferentially'' critical points.
In \cite{BDL} this fact is used to prove that such function has
only finitely many critical points of that sort. Both are direct
consequences of Theorem 1 (the latter even for all  definable
functions) if we define subdifferentially critical value  of a
function as just the critical value in the same sense as above.
Such a definition is justifiable as, according to the ``point''
subdifferential criterion for regularity (see e.g.
\cite{AI00,BM,RW}) a point is subdifferentially critical point of
a continuous function precisely when  the rate of surjection of
the function (considered as a mapping into $\R$) is zero at this
point. But we give an independent proof because, unlike the proof
of Theorem 1, it does not use the argument ad absurdum. Theorem 4
can also be used to get a direct proof of Theorem 1 for single
valued locally Lipschitz tame mappings using the fact that for
such mappings $\sur F(x)=0$ is equivalent to $\sd(y^*\circ
F)(x)=0$ for some $y^*$ with $\|y^*\|=1$.

\vskip 2mm

\thm 4. {\it Let $f$ be a definable  continuous function. Then $f$
is constant on every connected component of the set of its
critical points.}

\proof Let $u$ and $w$ be two different points belonging to the
same connected components of the set of critical points of $f$.
This set is definable since  so is the function $\sur f$ by
Proposition 1. Hence there is a definable curve $x(t),\ 0\le t\le
1$ joining $u$ and $w$ and lying completely in the set of critical
points of $f$. By definition for any $\ep>0$ and any $t\in[0,1]$
there is an $x$ and $\la >0$ such that $\| x-x(t)\|<\ep$,
$0<\la<\ep$ and $\Sur~ f(x)(\la)<\ep\la$. By the definable choice
theorem there are functions $z_{\ep}(t)$ and $\la_{\ep} (t)$
(definably depending on both variables) such that $\|
x(t)-z_{\ep}(t)\|<\ep$, \  $0<\la_{\ep}(t)<\ep$ and $\Sur~
f(z_{\ep}(t))(\la_{\ep}(t))<\ep\la_{\ep}(t)$ for all $t$ and all
$\ep$. If we set
$$
\mu_{\ep}^+(t)=\sup\{ f(z)-f(z_{\ep}(t)):\;
\|z-z_{\ep}(t)\|<\la_{\ep}(t)\},
$$
$$
 \mu_{\ep}^-(t)=\sup\{
f(z_{\ep}(t))-f(z):\; \|z-z_{\ep}(t)\|<\la_{\ep}(t)\},
$$
the latter amounts to
$$
\mu(t)=\min\{\mu_{\ep}^+(t),\mu_{\ep}^-(t)\}
<\ep\la_{\ep}(t).
\leqno (9)
$$

Applying Proposition 5 to each component of $x(t)-z_{\ep}(t)$, we
conclude that
$$
\int_0^1\| \dot x(t)-\dot z_{\ep}(t)\| dt\to 0. \leqno (10)
$$

The functions $z_{\ep}(\cdot)$ and $\la_{\ep}(\cdot)$ may have
points of discontinuity but by the uniform finiteness theorem the
number of such points is  bounded by the same constant for all
$\ep$. Note also that as $x(t)$ is continuous, the size of each of
jump of $z_{\ep}$ does not exceed $2\ep$.

We observe further that  $\mu_{\ep}^{\pm}$ and $\mu_{\ep}$ are
definable functions.  Therefore for any $\ep$ there are finitely
many points on $[0,1]$ such that between any pair of adjacent
points either $\mu_{\ep}(t)=\mu_{\ep}^+(t)$ or
$\mu_{\ep}(t)=\mu_{\ep}^-(t)$.

Thus there is a natural $N$ such that for any  $\ep>0$ there are
points  $\tau_i,\ i=1,...,k$,
$(0=\tau_0\le\tau_1<...<\tau_k\le\tau_{k+1}=1, \ k\le N)$ such
that on every interval $(\tau_i,\tau_{i+1})$ $z_{\ep}(t)$ and
$\la_{\ep(t)}$ are continuous and either
$\mu_{\ep}(t)=\mu_{\ep}^+(t)$ or $\mu_{\ep}(t)=\mu_{\ep}^-(t)$ for
all $t$ in the interval

As $f$ is continuous, the theorem will be proved if we show that
$$
|f(z_{\ep}(0))-f(z_{\ep}(1))|\to 0. \leqno (11)
$$
So fix an $\ep>0$ and let $\tau_i,\ i=1,...,k$,
$(0=\tau_0\le\tau_1<...<\tau_k\le\tau_{k+1}=1, \ k\le N)$ be the
points specified above. For any $i$ we set
$$
z_{\ep}(\tau_i)^+=\lim_{t\to\tau_i+0} z_{\ep}(t),\quad
z_{\ep}(\tau_i)^-=\lim_{t\to\tau_i-0} z_{\ep}(t).
$$

Take a certain interval $(\tau_i,\tau_{i+1})$ and assume for
instance that $\mu_{\ep}(t)=\mu_{\ep}^+(t)$ on this interval. This
means  for any $t\in (\tau_i,\tau_{i+1})$ there is a
$t'\in(t,\tau_{i+1})$ such that
$\|z_{\ep}(t')-z_{\ep}(t)\|<\la_{\ep}(t)$. Fix a $t$ and let $t^+$
be the upper bound of such $t'$. Then in the inequality above we
either get an equality at $t^+$ (as $z_{\ep}$ and $\la_{\ep}$ are
continuous on the interval) or $t^+$ coincides with the right end
of the interval. In the last case
$$
|f(z_{\ep}(\tau_{i+1})^-)-f(z_{\ep}(t))|\le\ep^2.
$$

If we get an equality at some $t^+\le \tau_{i+1}$, then, as
$\la_{\ep}(t)<\ep$, we have
$$
|f(z_{\ep}(t^+))-f(z_{\ep}(t))|\le\ep\la_{\ep}(t^+)=\ep\|
z_{\ep}(t^+)-z_{\ep}(t)\| \le\ep\int_{t}^{t^+}\|\dot
z_{\ep}(s)\|ds.
$$
We can ask about the upper bound $\tau^+$ of $t^+$ for which the
last inequality holds. The standard argument shows that either
this upper bound is $\tau_{i+1}$ or
$\|z_{\ep}(\tau_{i+1})-z_{\ep}(\tau^+)\|<\la_{\ep}(\tau^+)<\ep$.
Indeed, if the opposite inequality holds, then there is a $\tau\in
(\tau^+,\tau_{i+1}]$ such that
$\|z_{\ep}(\tau)-z_{\ep}(\tau^+)\|=\la_{\ep}(\tau^+)$ and
therefore
$$
|f(z_{\ep}(\tau)-f(z_{\ep}(\tau^+))|\le \ep
\int_{\tau^+}^{\tau}\|\dot z_{\ep}(s)\|ds
$$
and we arrive to a contradiction with the definition of $\tau^+$.
Thus we can conclude by stating that for any $t$ in the interval
$$
|f(z_{\ep}(\tau_{i+1})^-)-
f(z_{\ep}(t)|\le\ep(\int_t^{\tau_{i+1}}\|\dot z_{\ep}(s)\|ds
+\ep).
$$
and consequently, by continuity
$$
|f(z_{\ep}(\tau_{i+1})^-)-
f(z_{\ep}(\tau_i)^+|\le\ep(\int_{\tau_i}^{\tau_{i+1}}\|\dot
z_{\ep}(t)\|dt +\ep). \leqno (12)
$$
The same argument, with an obvious change, applies to intervals on
which $\mu_{\ep}(t)=\mu_{\ep}^-(t)$, and we can be sure that (12)
holds for each interval of the partition.

 Let finally $\omega (r)$ be the modulus of continuity of $f$
in a neighborhood of $x(\cdot)$. Then taking into account (12)
along with the fact that we have at most $N$ points of
discontinuity of $z_{\ep}$ and the jumps cannot exceed $2\ep$, we
find that
$$
|f(z_{\ep})(1))-f(z_{\ep}(0))|\le \ep(N+1)\big (\int_0^1\| \dot
x(t)\|dt +\ep + \int_0^1\| \dot x(t)-\dot z_{\ep}(t)\|dt\big )
+N\omega(2\ep)\big )
$$
from which (11) immediately follows in view of (10).
\qquad\qquad\qquad\qquad\qquad\black


\newpage

\begin{thebibliography}{99}

\bibitem{DAC}
D. d'Acunto, Valeurs critiques asymptotiques de fonctions
d\'efinissables dans un structure $o$-minimale, {\it Ann. Polon.
Math}, {\bf LXXV} (2000), 35-45.

\bibitem{SB}
S.M. Bates, Toward a precise smoothness hypothesis in Sard's
theorem, {\it Proc. Amer. Math. Soc.} {\bf 117} (1993), 279-283.

\bibitem{BDL}
J. Bolte, A. Daniilidis and A. Lewis, A nonsmooth Morse-Sard
theorem for subanalytic functions, preprint 2005.

\bibitem{BDL1}
J. Bolte, A. Daniilidis and A. Lewis,  The Lojasiewicz inequality
for nonsmooth subanalytic functions with applications to
subgradient dynamical systems, {\it SIAM J. Optimization}, to
appear

\bibitem{BDLS}
J. Bolte, A. Daniilidis and A. Lewis and M. Shiota, A Sard-type
theorem for Clarke critical values of subanalytic Lipschitz
continuous functions, preprint.

\bibitem{BCR}
J. Bochnak, M. Coste and M.-F. Roy, {\it Real Algebraic Geometry},
Springer 1998.

\bibitem{BZ}
J.M. Borwein and Q.J. Zhu, {\it Techniques of Variational
Analysis}, Springer 2005.

\bibitem{MC}
M. Coste, {\it An Introduction to $o$-Minimal Geometry}, Inst.
Rech. Math., Univ. de Rennes, 1999
(http://name.math.univ-rennes1.fr/michel.coste/polyens/OMIN.pdf)

\bibitem{DMT}
E. De Giorgi, A. Marino and M. Tosques, Problemi di evoluzione in
spazi metrici e curve di massima pendenza, {\it Atti Acad. Nat.
Lincei, Rend. Cl. Sci. Fiz. Mat. Natur.} {\bf 68} (1980), 180-187.

\bibitem{LVD}
L. van den Dries, {\it Tame Topology and O-minimal Structures},
Cambrifge Univ. Press 1998.

\bibitem{AG}
A. Grothendieck, \ Sketch of a proposal, in {\it Geometric Galois
Actions}, L. Schneps and P. Lochak, eds. Cambridge Univ. Press
1997.

\bibitem{AI00}
A.D. Ioffe, Metric regularity and subdifferential calculus, {\it
Uspekhi Matem. Nauk} {\bf 55:3} (2000), 103-162, English
translation: {\it Russian Math. Surveys} {\bf 55} (2000), 501-558.

\bibitem{KOS}
 K. Kurdyka, P. Orro and S. Simon, Semialgebraic Sard
theorem for generalized critical values, {\it J. Differential
Geometry} {\bf 56} (2000), 67-92.

\bibitem{AL}
A.S. Lewis, Active sets, nonsmoothness and sensitivity, {\it SIAM
J. Optimization} {\bf 13} (2003, 702-725.


\bibitem{RM}
R. Mifflin, Semismooth and semiconvex functions in constrained
optimization, {\it Math. Operation res.} {\bf 2} (1977), 191-207.

\bibitem{BM}
B.S. Mordukhovich, {\it Variational Analysis and Generalized
Differentiation}, vol 1. Springer 2005.

\bibitem{AN}
A. Norton, Functions not constant on fractal quasi-arcs of
critical points, {\it Proc. Amer. Math. Soc} {\bf 106} (1989),
397-405.

\bibitem{LR}
L.  Rifford, A Morse-Sard theorem for the distance function on
Riemannian manifolds, {\it Manuscripta Math.} {\bf 113} (2004),
251-265.

\bibitem{RT}
R.T. Rockafellar, Favorable classes of Lipschitz continuous
functions in subgradient optimization, in {\it Progress in
Non-Differentiable Optimization}, E. Nurminski, editor, Pergamon
Press 1981.

\bibitem{RW}
R.T. Rockafellar and R.J.B. Wets, {\it Variational Analysis},
Springer 1998

\bibitem{AR}
A. Rohde, On Sard's theorem for nonsmooth functions {\it Numer.
Funct. Anal. Optim.} {\bf 9-10} (1997), 1023-1039.

\bibitem{BT}
B. Teissier, Tame and stratified objects, in {\it Geometric Galois
Actions}, L. Schneps and P. Lochak, eds. Cambridge Univ. Press
1997, pp. 231-242.

\bibitem{HW}
H. Whitney, A function not constant on a connected set of critical
points, {\it Duke math. J.} {\bf 1} (1935), 514-517.

\bibitem{YY1}
Y. Yomdin,  Maxima of smooth families III: Morse-Sard theorem,
preprint, MOI, Bonn 1984


\bibitem{YY3}
Y. Yomdin, Surjective mapping whose differential is nowhere
surjective, {\it Proc. Amer. Math. Soc} {\bf 111} (1991),267-270.



\end{thebibliography}
\end{document}